\documentclass[12pt]{article}
\usepackage{amsfonts}
\usepackage{amssymb}
\usepackage{CJK}
\usepackage{mathrsfs}%有更多的数学符号，比如花体 $\mathscr P$

\usepackage{amsfonts,amssymb, amsmath, latexsym}
\usepackage{CJK}
\usepackage{pst-poly}     % From pstricks/contrib/pst-poly
\usepackage[left,pagewise,mathlines]{lineno}
\usepackage{cite}

\newtheorem{thm}{Theorem}[section]

\newtheorem{lem}[thm]{Lemma}

\newenvironment{pf}[1][Proof]{\noindent\textbf{#1.} }{\hfill\rule{1mm}{2mm}}

\begin{document}

\title{Many-to-many disjoint paths in hypercubes with faulty vertices \thanks{The work was supported
by NNSF of China (Nos 10711233, 11101378) and ZSDZZZZXK08.}}

\author
{{\large Xiang-Jun Li$^a$ \quad Bin Liu$^a$ \quad Meijie Ma$^b$\quad
Jun-Ming Xu$^a$\footnote{Corresponding
author: xujm@ustc.edu.cn (J.-M. Xu)}} \\
{\small $^a$School of Mathematical Sciences}\\
{\small  University of Science and Technology of China}\\
{\small Wentsun Wu Key Laboratory of CAS}\\
{\small Hefei 230026, China}\\
{\small $^b$Department of Mathematics}\\
{\small Zhejiang Normal University}\\
{\small Jinhua 321004, China}  \\
 }

\date{}

\maketitle

\begin{abstract}
This paper considers the problem of many-to-many disjoint paths in
the hypercube $Q_n$ with $f$ faulty vertices and obtains the
following result. For any integer $k$ with $1\leq k\leq n-2$, any
two sets $S$ and $T$ of $k$ fault-free vertices in different parts
of $Q_n\ (n\geq 3)$, if $f\leq 2n-2k-3$ and each fault-free vertex
has at least two fault-free neighbors, then there exist $k$ fully
disjoint fault-free paths linking $S$ and $T$ which contain at least
$2^n-2f$ vertices. This result improves some known results in a
sense.
\end{abstract}



%%%%%%%%%%%%%%%%%%%%%%%%%%%%%%%%%%%%%%%%%%%%%%%%%%%%%%%%%%%%%%%%%%%%%%%%%%%%%%%%%%%%%%%%%%%%%%%%%%%%%%%%%%%%%%%%%%
\section{Induction}

The $n$-dimensional hypercube $Q_n$ is a graph whose vertex-set
consists of all binary vectors of length $n$, with two vertices
being adjacent whenever the corresponding vectors differ in exactly
one coordinate. There is a large amount of literature on
graph-theoretic properties of hypercubes (e.g., see the
comprehensive survey paper on early results~\cite{hhw88}) and recent
results~\cite{x01} as well as on their applications in parallel
computing (e.g., see~\cite{l92}).

One of most central issues in various high performance communication
networks or parallel computing systems is to find a cycle or a path
of given length in $Q_n$ (see~\cite{dk09, fg09, gs09} and the survey
paper~\cite{xm09}). To find node-disjoint paths concerned with
fault-tolerant routings among nodes has received much research
attention recently (see, for example, \cite{c09, m10, pkl09} and
references cited therein). In this paper, we consider the problem of
node disjoint paths of $Q_n$.

It is well known that there are $n$ internally disjoint paths
between any two vertices $u$ and $v$ in $Q_n$, but they do not
always contain all vertices. In 2004, Chang {\it et
al.}~\cite{clhh04} showed that for any integer $k$ with $1\leq k\le
n$ and any two vertices $u, v$ from different partite sets of $Q_n$,
there are $k$ internally disjoint $uv$-paths in $Q_n$ which contains
all vertices. In 2007, Chang {\it et al.}~\cite{cshh07} further
showed that for any integer $k$ with $1\leq k\le n-4$ and any two
vertices $u, v$ from different partite sets of $Q_n\ (n\geq 5)$,
there are $k$ internally disjoint $uv$-paths in $Q_n$ which contain
all vertices and the difference of lengths of any two paths is most
two. In 2009, Chen~\cite{c09} showed that for any integer $k$ with
$1\leq k\leq n-1$, any two sets $S$ and $T$ of $k$ fault-free
vertices in different partite sets in $Q_n$ with $f$ faulty vertices
and $h$ faulty edges, if $f+h\leq n-k-1$, then there exist $k$
disjoint fault-free paths linking $S$ and $T$ in $Q_n$ which contain
at least $2^n-2f$ vertices. In this paper, we obtain the following
results.

\noindent{\bf Theorem}\ {\it Let $Q_n\ (n\geq 3)$ be an
$n$-dimensional hypercube with $f$ faulty vertices, and let $k$ be
an integer with $1\leq k\leq n-2$. If $f\le 2n-2k-3$ and each
fault-free vertex has at least two fault-free neighbors then, for
any two sets $S$ and $T$ of $k$ fault-free vertices in different
partite sets, there exist $k$ disjoint fault-free paths linking $S$
and $T$ in $Q_n$ which contain at least $2^n-2f$ vertices.}

Clearly, our theorem partly improves the above-mentioned results.
The proof of our result is in Section 3. Section 2 gives some
notations and lemmas.

\section{Notations and Lemmas}

For graph-theoretical terminology and notation not defined here, we
follow \cite{x03}. Let $G=(V,E)$ be a connected simple graph, where
$V=V(G)$ is the vertex-set and $E=E(G)$ is the edge-set of $G$. For
$uv\in E(G)$, we call $u$ (resp. $v$) is a neighbor of $v$ (resp.
$u$). A $uv$-path is a sequence of adjacent vertices, written as
$\left<v_0,v_1,v_2,\cdots,v_m\right
>$, in which $u=v_0, v=v_m$ and all the vertices
$v_0,v_1,v_2,\cdots,v_m$ are different from each other. The length
of a path $P$ is the number of edges in $P$. Let $d_G(u,v)$ be the
length of a shortest $uv$-path in $G$, called the distance between
$u$ and $v$ in $G$.
%For any two subsets $S$ and $T$ of $V(G)$, the distance between
%$S$ and $T$ is defined as $d_G(S,T)=\min\{d_G(x,y):\ x\in S, y\in T\}$.
For a path $P=\left<v_0,v_1,\cdots,v_i, v_{i+1},
\cdots, v_m\right>$, we can write
$P=P(v_0,v_i)+v_iv_{i+1}+P(v_{i+1},v_m)$, and the notation
$P-v_iv_{i+1}$ denotes the subgraph obtained from $P$ by deleting
the edge $v_iv_{i+1}$. Two paths are disjoint if they have no
vertices in common. Given two disjoint sets $S$ and $T$ of $k$
vertices, if there exist $k$ disjoint paths between $S$ and $T$, we
call these paths to be $k$ disjoint $ST$-paths.

The $n$-dimensional hypercube $Q_n$ is a graph with $2^n$ vertices,
each vertex denoted by an $n$-bit binary string $u=u_nu_{n-1}\cdots
u_2u_1$. Two vertices are adjacent if and only if their strings
differ in exactly one bit position. It has been proved that $Q_n$ is
a vertex- and edge-transitive bipartite graph.

By definition, for any $k\in\{1,2,\cdots,n\}$, $Q_n$ can be
expressed as $Q_n=L_k\odot R_k$, where $L_k$ and $R_k$ are two
subgraphs of $Q_n$ induced by the vertices with the $k$ bit position
is 0 and $1$, respectively, which are isomorphic to $Q_{n-1}$,
linking by an edge a vertex in $L_k$ to a vertex in $R_k$ if they
differ in only the $k$-th bit. Without loss of generality, we write
$Q_n=L\odot R$. For convenience, for a vertex $u$ in $L$, we use
$u_R$ to denote its only neighbor in $R$. Similarly, for a vertex
$v$ in $R$, we use $v_L$ to denote its only neighbor in $L$.
Clearly, for any two vertices $u_L$ and $v_L$ in $L$,
$d_L(u_L,v_L)=d_R(u_R,v_R)$.

Let $F$ denote a set of faulty vertices in $Q_n$, and $f=|F|$. When
$Q_n=L\odot R$, we denote $f_L=|F\cap L|$ and $f_R=|F\cap R|$. A
subgraph of $Q_n$ is {\it fault-free} if it contains no vertices in
$F$. For two subsets $A,B$ of $V(Q_n)$, let
$d(A,B)=\min\{d(x,y):x\in A, y \in B\}$.

\begin{lem}\label{lem2.1} \textnormal{(Kueng {\it et al.}~\cite{klht09},
2009)} If $f\leq 2n-5$ and each fault-free vertex of $Q_n (n\geq 3)$
has at least two fault-free neighbors then, for any two distinct
fault-free vertices $x$ and $y$ with distance $d$, there is
a fault-free $xy$-path containing at least $2^{n}-2f$ vertices if
$d$ is odd and $2^{n}-2f-1$ vertices if $d$ is even.
\end{lem}

\begin{lem}\label{lem2.2} \textnormal{(Chen~\cite{c09}, 2009)}
For any integer $k$ with $1\leq k\leq n-1$, if $f\leq n-k-1$ then,
for any two sets $S$ and $T$ of $k$ fault-free vertices in different
partite sets in $Q_n\ (n\geq 2)$, there exist $k$ disjoint fault-free
$ST$-paths which contain at least $2^n-2f$ vertices.
\end{lem}

\begin{lem}\label{lem2.3} \textnormal{(Dvo\v{r}\'{a}k~\cite{d05}, 2005)}
Let $(x, y)$ and $(u, v)$ be two disjoint pairs of vertices with odd
distance in $Q_n$. If $x$ and $y$ are adjacent in $Q_n$ with $n\ge
3$, then there exists a $uv$-path containing all vertices in
$Q_n-\{x,y\}$ unless $n=3, d(u,v)=1$ and $d(\{x,y\},\{u,v\})=2$.
%
%Let $(x, y)$ and $(u, v)$ be two disjoint pairs of vertices
%with odd distance in $Q_n\, (n\ge2)$. Then there exist two disjoint
%$xy$-path $P_1$ and $uv$-path $P_2$ containing all vertices in $Q_n$.
%such that $V(P_1)\cap V(P_2)=\emptyset, V(P_1)\cup V(P_2)=V(Q_n)$.
%Moreover, if $d(x,y)=1$, then $P_1$ can be chosen such that $P_1$ is
%exactly the edge $xy$, unless $n=3$, $d(u,v)=1$ and
%$d(\{x,y\},\{u,v\})=2$.
\end{lem}

%Simmons [10]
%introduced the concept of hamiltonian laceability on bipartite
%graphs due to its importance. A bipartite graph is defined to
%be hamiltonian laceable if (a) it is equitable and whenever x
%and y are two vertices from different partite sets, there exists
%an x,y-hamiltonian path; or else (b) it is nearly equitable and
%whenever x and y are two vertices in the larger partite set,
%there exists an x,y-hamiltonian path.

%[10] G. Simmons, Almost all n-dimensional rectangular lattices are hamilton
%laceable, Congr. Numer., 21 (1978) 649-661

%Lemma 3: [12] For $n\geq 4$, $Q_n-\{x,y\}$ is hamiltonian
%laceable, where $x$ and $y$ are any two vertices from different
%partite sets.

%[12] C. M. Sun, C. K. Lin, H. M. Huang, L. H. Hsu, Mutually independent
%hamiltonian paths and cycles in hypercubes, J. Interconnection
%Networks, 7 (2) (2006) 235-255.

%\begin{lem}\label{lem2.3a} \textnormal{(Sun {\it et al.}~\cite{slhh06}, 2006)}
%Let $(x, y)$ and $(u, v)$ be any two disjoint pairs of vertices
%with odd distance in $Q_n\, (n\ge4)$. Then there exists a $uv$-path containing all vertices
%in $Q_n-\{x,y\}$.
%\end{lem}

\begin{lem}\label{lem2.4}
Let $x$ and $y$ be two adjacent vertices of $Q_n\,(n\ge 4)$. Then,
for any integer $k$ with $1\leq k\leq n-2$ and for any two sets $S$
and $T$ of $k$ vertices in different partite sets in $Q_n-\{x,y\}$,
there exist $k$ disjoint $ST$-paths containing all vertices in
$Q_n-\{x, y\}$.
 \end{lem}

\begin{pf}
By induction on $k$. If $k=1$, the result is true by
Lemma~\ref{lem2.3}. We assume that the result holds for any integer
fewer than $k$, and consider the case of $k\,(\geq 2)$. Let
$\{X,Y\}$ be a bipartition of $Q_n$ and $Q_n=L\odot R$ such that the
edge $xy\in E(L)$, $S=\{s_1,s_2,\ldots, s_k\}$ and $T=\{t_1,
t_2,\ldots,t_k\}$ be two sets of $k$ vertices in different partite
sets in $Q_n-\{x,y\}$, and let
  $$
  \begin{array}{rl}
  &S_L= S\cap L,\ T_L=T\cap L\ \ {\rm and}\ \ p=|S_L|,\ q=|T_L|,\\
  &S_R= S\cap R,\ T_R=T\cap R.
  \end{array}
  $$
Without loss of generality, we can assume $S\cup\{x\}\subseteq X$,
$T\cup\{y\}\subseteq Y$, $p\ge q$ and $0\le q \le k-1$. If $p > q$,
let
 $$
 \begin{array}{rl}
 &U_L=\{u_1,\ldots, u_{p-q}\}\subseteq Y\cap V(L-(T_L\cup\{y\}),\\
 &U_R=\{v_1,\ldots, v_{p-q}\}\subseteq X\cap V(R)\ {\rm with}\ u_iv_i\in E(Q_n)\ {\rm for}\ 1\le i\le p-q.
 \end{array}
 $$ If
$p=q$, let $U_L=U_R=\emptyset$.

The process of our induction steps stronly depends on
Lemma~\ref{lem2.3}. When $n=4$, $L\cong R\cong Q_3$, while in $Q_3$,
we have to consider some exceptional cases when we apply
Lemma~\ref{lem2.3}. Thus, we first prove the case of $n\ge 5$ by
considering the following three cases.

{\bf Case 1}. $p=0$. In this case, both $S$ and $T$ are in $R$. By
Lemma~\ref{lem2.2}, there exist $k$ disjoint $ST$-paths
$P_1,P_2,\ldots,P_k$ containing all vertices in $R$. We can choose
an edge $u_Lv_L$ in $L-\{x,y\}$ such that $u_Rv_R$ is an edge in
some path $P_i$. Since $n\ge 5$, by Lemma~\ref{lem2.3}, there exists
a $u_Lv_L$-path $P_0$ containing all vertices in $L-\{x,y\}$. Let
 $$
 P'_i=P_i-u_Rv_R+u_Ru_L+v_Lv_R+P_0.
 $$
Replacing $P_i$ by $P'_i$ obtains $k$ disjoint $ST$-paths
$P_1,\ldots,P_{i-1},P'_i,P_{i+1},\ldots,P_k$ as required.

{\bf Case 2}. $1\le p\le k-1$. Let $T'_L=T_L\cup U_L$ and
$S'_R=S_R\cup U_R$. Then $|T'_L|=q+(p-q)=p=|S_L|$ and
$|S'_R|=(k-p)+(p-q)=k-q=|T_R|$. Since $n\ge 5$, by the induction
hypothesis, there are $p$ disjoint $S_LT'_L$-paths containing all
vertices in $L-\{x,y\}$, in which let $P_1,P_2,\ldots,P_{q}$ be
$S_LT_L$-paths, and $P'_1,P'_2,\ldots,P'_{p-q}$ be $S_LU_L$-paths,
without loss of generality, say $P'_i$ connecting one vertex in
$S_L$ to $u_i$. By Lemma~\ref{lem2.2}, there are $k-q$ disjoint
$S'_RT_R$-paths containing all vertices in $R$, in which let
$P''_1,P''_2,\ldots,P''_{p-q}$ be $U_RT_R$-paths, say $P''_i$
connecting $v_i$ to one vertex in $T_R$, and
$P_{p+1},P_{p+2},\ldots,P_{k}$ be $S_RT_R$-paths. Let
 $$
 P_{q+i}=P'_i+u_iv_i+P''_i\ \ {\rm for\ each}\ i=1,2,\ldots,p-q.
 $$
Then $P_1,P_2,\ldots,P_k$ are $k$ disjoint $ST$-paths containing all
vertices in $Q_n-\{x, y\}$.

{\bf Case 3}. $p=k$.  By the hypothesis of $q \le k-1$,
$T_R\ne\emptyset$, say $t_k\in T_R$. Let $s_k\in S$, $S'=S- \{s_k\},
U'=U_L- \{u_{k-q}\}$ and $T'=T_L\cup U'$. Since $n\ge 5$, by the
induction hypothesis, there exist $k-1$ disjoint $S'T'$-paths
containing all vertices in $L-\{x,y\}$, in which, let
$P_1,P_2,\ldots, P_q$ be $S'T_L$-paths, and
$P'_1,P'_2,\ldots,P'_{k-1-q}$ be $S'U'$-paths, where $P'_{i}$
connects  one vertex of $S'$ to $u_i$ for each $i=1,2,\ldots,k-1-q$.

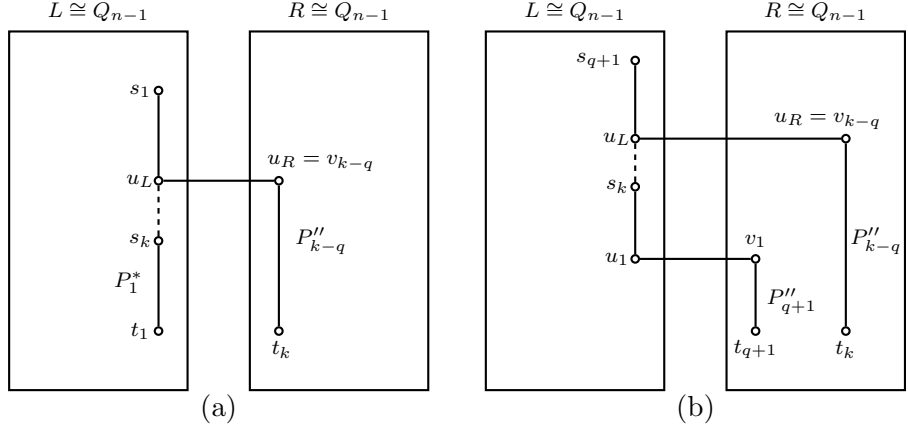
\begin{figure}[h]
\begin{center}

\begin{pspicture}(0,1)(5.2,6)
\psset{unit=0.8, radius=.08}

\psframe(0,1)(3,7) \psframe(4,1)(7,7)

\Cnode(2.5,6){s1}\rput(2.2,6){\scriptsize$s_1$}

\Cnode(2.5,2){t1} \rput(2.2,2){\scriptsize$t_1$}
\rput(2,2.8){\scriptsize$P_1^*$}

\Cnode(2.5,4.5){uL}
\rput(2.2,4.5){\scriptsize$u_L$}\Cnode(2.5,3.5){vL}
\rput(2.2,3.5){\scriptsize$s_k$} \ncline{s1}{uL}\ncline{t1}{vL}
\ncline[linestyle=dashed, dash=2pt 2pt]{uL}{vL}

\Cnode(4.5,4.5){uR}
\rput(5.2,4.8){\scriptsize$u_R=v_{k-q}$}\Cnode(4.5,2){vR}
\rput(4.55,1.7){\scriptsize$t_k$}

\ncline{uR}{vR} \ncline{uR}{uL}

\rput(5.2,3.5){\scriptsize$P_{k-q}''$}

\rput(1.5,7.3){\scriptsize$L\cong
Q_{n-1}$}\rput(5.5,7.3){\scriptsize$R\cong Q_{n-1}$}
\rput(3.5,0.7){\small (a)}
\end{pspicture}
\begin{pspicture}(-1,1)(4.2,6)
\psset{unit=0.8, radius=.08}

\psframe(0,1)(3,7) \psframe(4,1)(7,7)

\Cnode(2.5,6.5){s1}\rput(1.9,6.5){\scriptsize$s_{q+1}$}

\Cnode(2.5,3.2){t1} \rput(2.2,3.2){\scriptsize$u_1$}
%\rput(2,2.8){\scriptsize$P_1^*$}

\Cnode(2.5,5.2){uL}
\rput(2.2,5.2){\scriptsize$u_L$}\Cnode(2.5,4.4){vL}
\rput(2.2,4.4){\scriptsize$s_k$} \ncline{s1}{uL}\ncline{t1}{vL}
\ncline[linestyle=dashed, dash=2pt 2pt]{uL}{vL}

\Cnode(6,5.2){uR}
\rput(5.7,5.5){\scriptsize$u_R=v_{k-q}$}\Cnode(6,2){vR}
\rput(6,1.7){\scriptsize$t_k$}

\Cnode(4.5,2){tq} \rput(4.5,1.7){\scriptsize$t_{q+1}$}
\Cnode(4.5,3.2){v1} \rput(4.5,3.5){\scriptsize$v_1$}

\ncline{tq}{v1}\ncline{v1}{t1}

\ncline{uR}{vR} \ncline{uR}{uL}

\rput(6.5,3.5){\scriptsize$P_{k-q}''$}
\rput(5.1,2.5){\scriptsize$P_{q+1}''$}

\rput(1.5,7.3){\scriptsize$L\cong
Q_{n-1}$}\rput(5.5,7.3){\scriptsize$R\cong Q_{n-1}$}
\rput(3.5,0.7){\small (b)}
\end{pspicture}

\end{center}
\caption{\label{f4}\footnotesize {Illustrations for the proof of
Case 3 of Lemma~\ref{lem2.4}.}}
\end{figure}

Assume that $s_k$ is in some path connecting $s_i$ to some vertex
$t'$ in $T'$. Let $u_L$ be the neighbor of $s_k$ in the path closer
$s_i$ and $u_R$ be the neighbor of $u_L$ in $R$, and let
$u_R=v_{k-q}$. By Lemma~\ref{lem2.2}, there are $k-q$ disjoint
$U_RT_R$-paths $P''_1,P''_2,\ldots,P''_{k-1-q},P''_{k-q}$ that
contain all vertices in $R$. Without loss of generality, assume that
$P''_{i}$ connects $v_i$ and $t_{q+i}$ for each $i=1,2,\ldots, k-q$.
Let
 $$
 P_{q+i}=P'_i+u_iv_i+P''_i\ \ {\rm for\ each}\ i=1,2,\ldots,k-1-q.
 $$

If $t'\in T_L$, say $t'=t_1$ and $s_k$ in $P_1$ connecting $s_1$ to
$t_1$ (see Fig.~\ref{f4} (a)), let
 $$
 \begin{array}{rl}
  &P^*_1=P_1(s_k,t_1),\\
 &P_k=P_1(s_1,u_L)+u_Lu_R+P''_{k-q}.
 \end{array}
 $$
Then $P^*_1,P_2,\ldots, P_{k-1},P_k$ are $k$ disjoint $ST$-paths
containing all vertices in $Q_n-\{x, y\}$.

If $t'\in U'$, say $t'=u_1$ and $s_k$ in $P'_{q+1}$ connecting
$s_{q+1}$ to $u_1$ (see Fig.~\ref{f4} (b)), let
$$
 \begin{array}{rl}
  &P^*_{q+1}=P'_{q+1}(s_k,u_1)+u_1v_1+P''_{q+1},\\
 &P_k=P'_{q+1}(s_1,u_L)+u_Lu_R+P''_{k-q}.
 \end{array}
 $$
Then $P_1,\ldots, P_q, P^*_{q+1}, P_{q+2}, \ldots, P_{k-1},P_k$ are
$k$ disjoint $ST$-paths containing all vertices in $Q_n-\{x, y\}$.

Summing up the three cases, we prove the lemma when $n\ge 5$.

When $n=4$, we have $k=2$. In this case, if there exists $i\,(1\le
i\le 2)$ such that $ d(\{s_i,t_i\},\{x,y\})=2$, we can choose
$Q_4=L\odot R$ such that $xy\in E(L)$ and $\{s_i,t_i\}\subseteq R$.
Then, by Lemma~\ref{lem2.3}, we can direct verify that the lemma
holds, and the details are omitted.
\end{pf}

\section{Proof of Theorem}

The proof proceeds by induction on $n\geq 3$. If $k=1$, then $f\leq
2n-5$. The theorem follows by Lemma~\ref{lem2.1}. If $k=n-2$, then
$f\leq 1$. The theorem follows by Lemma~\ref{lem2.2}. Thus, the
theorem holds for $3\leq n\leq 4$. In the following discussion, we
assume $n\geq 5$ and $2\leq k\le n-3$.

Let $Q_n$ be an $n$-dimensional hypercube, $\{X,Y\}$ be a
bipartition of $V(Q_n)$, $S=\{s_1,s_2,\ldots, s_k\}\subset X$ and
$T=\{t_1, t_2,\ldots,t_k\}\subset Y$ be any two sets of $k$
fault-free vertices. Our aim is to construct $k$ disjoint fault-free
$ST$-paths containing at least $2^n-2f$ vertices.

Since $k\geq 2$, we have $f\leq 2n-7$. It is easy to see that there
is some $j\in \{1,2,\ldots,n\}$ such that each fault-free vertex in
$L_{j}$ and $R_{j}$ has at least two fault-free neighbors. Let
$Q_n=L\odot R$, where $L=L_{j}$ and $R=R_{j}$. Let
  $$
  \begin{array}{rl}
  &S_L= S\cap L,\ T_L=T\cap L\ \ {\rm and}\ \ p=|S_L|,\ q=|T_L|,\\
  &S_R= S\cap R,\ T_R=T\cap R.
  \end{array}
  $$
By symmetry, we may assume $p\geq q$.

We construct $k$ required $ST$-paths by considering three cases.

\vskip6pt {\bf Case 1.}\quad  $q=k$ or $p=0$.

In this case, two sets $S$ and $T$ both are in $L$ or $R$. By
symmetry, we consider the former one, that is, $S$ and $T$ both are
in $L$. There are two subcases.

\vskip6pt

{\it Subcase 1.1.}\quad $f_L\leq 2n-2k-5$.

Since $f_L\leq 2n-2k-5=2(n-1)-2k-3$ and $k\leq n-3=(n-1)-2$, by the
induction hypothesis, in $L$ there are $k$ disjoint $ST$-paths
 \begin{equation}\label{e3.1}
 P_1',P_2',\ldots,P_k'
 \end{equation}
containing at least $2^{n-1}-2f_L$ vertices. Note that, when $n\geq
5$ and $f\leq 2n-7$,
 $$
 \begin{array}{rl}
 \sum\limits_{i=1}^k|E(P_i')|&=\sum\limits_{i=1}^k|V(P_i')|-k\\
 &\geq  2^{n-1}-2f_L-k\\
 &\geq 2^{n-1}-2f-(n-3)\\
 &>4(2n-7)-2f\\
 &\geq 2f.
 \end{array}
 $$
Therefore, there is an edge $u_Lv_L$ in some path, say $P_1'$, such
that $u_R$ and $v_R$ both are fault-free. Since $f_R\leq f\leq
2(n-1)-5$, by Lemma~\ref{lem2.1}, in $R$ there is an $u_Rv_R$-path
$P_R$ containing at least $2^{n-1}-2f_R$ vertices. Let
 \begin{equation}\label{e3.2}
 P_1=P_1'-u_Lv_L+u_Lu_R+P_R+v_Rv_L,\ \ P_i=P_i'\ \ {\rm for\ each}\ i=2,3,\dots, k.
 \end{equation}
Then $P_1,P_2,\ldots,P_k$ are $k$ disjoint $ST$-paths containing at
least $2^n-2f$ vertices in $Q_n$.

\vskip6pt {\it Subcase 1.2.}\quad $f_L\geq 2n-2k-4$.

In this case, $f_R\leq 1$ since $f\leq 2n-2k-3$. We consider two
subcases according as whether $f_R$ is equal to 1 or not,
respectively.

{\it Subcase 1.2a.}\quad $f_R=1$.

In this case, $f=2n-2k-3$ and $f_L=2n-2k-4$. All paths in
(\ref{e3.1}) contain at most one faulty vertex. If they contain no
faulty vertex, then paths defined in (\ref{e3.2}) are as required.
Assume that some path, say $P_1'$, contains a faulty vertex $w$.
Without loss of generality, $P_1'$ connects $s_1$ and $t_1$.
%Then $P_1'$ contains at least $4$ vertices.
Note that, in this case,
paths in (\ref{e3.1}) contain at least $2^{n-1}-2(f_L-1)$ vertices
in $L$. Let $u_L$ and $v_L$ be two neighbors of $w$ in $P_1'$, where
$u_L$ is in $P_1'(s_1,w)$ and $v_L$ is in $P_1'(w,t_1)$. Then the
distance of $u_L$ and $v_L$ is even.

If $u_R$ and $v_R$ both are fault-free then, by Lemma~\ref{lem2.1},
there is a fault-free $u_Rv_R$-path $P_R$ containing at least
$2^{n-1}-2f_R-1$ vertices in $R$. Let
 $$
 P_1^*=P_1'(s_1,u_L)+u_Lu_R+P_R+ v_Rv_L+P_1'(v_L,t_1).
 $$
Then, replacing $P_1$ in (\ref{e3.2}) by $P_1^*$ yields $k$ disjoint
fault-free $ST$-paths with at least $2^n-2f$ vertices in $Q_n$.

If one of $u_R$ and $v_R$ is faulty vertex, say $u_R$, then $u_R$ is
the only faulty vertex in $R$. If $u_L\ne s_1$, let $z_L$ be the
fault-free neighbor of $u_L$ in $P_1'(s_1,u_L)$ (maybe $z_L=s_1$),
then $z_R$ is fault-free in $R$, and the distance $z_R$ and $v_R$ is
odd. By Lemma~\ref{lem2.1}, there is a fault-free $z_Rv_R$-path
$P_R$ containing at least $2^{n-1}-2f_R$ vertices in $R$. Let
 $$
 P_1^*=P_1'(s_1,z_L)+z_Lz_R+P_R+ v_Rv_L+P_1'(v_L,t_1).
 $$
Replacing $P_1$ in (\ref{e3.2}) by $P_1^*$ yields $k$ disjoint
fault-free $ST$-paths with vertices at least
 $$
 (2^{n-1}-2(f_L-1)-|\{w,u_L\}|)+(2^{n-1}-2f_R)=2^n-2f.
 $$

We now assume $u_L=s_1$ and $u_R$ is the only faulty vertex in $R$.
Note that $u_L$ has at least two fault-free neighbor in $L$.

If there is a fault-free neighbor $z_L$ of $u_L$ not in $P_i'$ for
each $i=1,2,\ldots,k$, we consider two vertices $z_R$ and $v_R$ in
$R$. They are free-fault in $R$ since $u_R$ is the only faulty
vertex in $R$. Since $d(z_L,v_L)=3$, the distance
between $z_R$ and $v_R$ is odd. Let $P_R$ be a $z_Rv_R$-path with
at least $2^{n-1}-2f_R$ vertices in $R$, and let
 $$
 P_1^*=s_1z_L+z_Lz_R+P_R+v_Rv_L+P_1'(v_L,t_1).
 $$
Replacing $P_1$ in (\ref{e3.2}) by $P_1^*$ yields $k$ disjoint
fault-free $ST$-paths with vertices at least
 $$
 (2^{n-1}-2(f_L-1)-|\{w\}|)+(2^{n-1}-2f_R)>2^n-2f.
 $$

If all the fault-free neighbors of $u_L$ are in $P_1'$, then one of
them, say $z_L$, is not $t_1$. Clearly, $z_L\ne v_L$ since $Q_n$
contains no triangles. Let $z_L'$ be the neighbor of $z_L$ in
$P_1'(z_L,t_1)$ (maybe $z_L'=t_1$), then $z_R'$ is fault-free in
$R$. Since both $v_L$ and $z_L'$ are in $X$, the distance between
$v_R$ and $z_R'$ is even. Let $P_R$ be a $v_Rz_R'$-path with at
least $2^{n-1}-2f_R-1$ vertices in $R$, and let
 $$
 P_1^*=s_1z_L+P_1'(z_L,v_L)+v_Lv_R+P_R+z_R'z_L'+P_1'(z_L',t_1).
 $$
Replacing $P_1$ in (\ref{e3.2}) by $P_1^*$ yields $k$ disjoint
fault-free $ST$-paths with vertices at least
 $$
 (2^{n-1}-2(f_L-1)-|\{w\}|)+(2^{n-1}-2f_R-1)=2^n-2f.
 $$

\begin{figure}[h]
\begin{center}

\begin{pspicture}(0,1)(6.5,7.5)
\psset{radius=.06}

\psframe(0,1)(3,7) \psframe(4,1)(7,7)
%
%\psframe[framearc=0.8](0.3,1.5)(2.8,2.4)
%\psframe[framearc=0.8](0.3,5.6)(2.8,6.5)

\psframe(0,1)(3,7) \psframe(4,1)(7,7)

\rput(1.5,7.3){\scriptsize$L\cong
Q_{n-1}$}\rput(5.5,7.3){\scriptsize$R\cong Q_{n-1}$}

\Cnode(4.5,5.2){zR}
\rput(4.5,4.9){\scriptsize$z_R'$}\Cnode(4.5,3){vR}
\rput(4.5,3.2){\scriptsize$v_R$}

\rput(5.3,4.2){\scriptsize$P_R$}

\nccurve[angleA=50,angleB=-50]{zR}{vR} \ncline{uR}{uL}

\Cnode(1.2,6.4){s2} \rput(0.9,6.4){\scriptsize$s_2$}
\Cnode(1.2,5.2){zL} \rput(1.45,5.45){\scriptsize$z_L'$}
\Cnode(1.2,4.5){zl} \rput(1.45,4.3 ){\scriptsize$z_L$}
\Cnode(2.5,4.5){s1} \rput(2.76,4.5 ){\scriptsize$s_1$}
\Cnode(2.5,3){vL} \rput(2.76,3.2){\scriptsize$v_L$}

\rput(0.9,4){\scriptsize$P_2'$} \rput(2.2,2.8){\scriptsize$P_1'$}

\Cnode(2.5,1.8){t1} \rput(2.5,1.5){\scriptsize$t_1$}
\Cnode(1.2,1.8){t2} \rput(1.2,1.5){\scriptsize$t_2$}

\Cnode[linestyle=dashed, dash=1pt 1pt](2.5,3.75){w} \rput(2.76,3.75
){\scriptsize$w$}

\ncline{t1}{vL}\ncline{t2}{zl}

\ncline{vR}{vL} \ncline{s2}{zL}\ncline{s1}{zl}\ncline{zR}{zL}

\ncline[linestyle=dashed, dash=2pt 2pt]{w}{s1}
\ncline[linestyle=dashed, dash=2pt 2pt]{w}{vL}

\ncline[linestyle=dashed, dash=2pt 2pt]{zl}{zL}

%\rput(3.5,0.7){\small (b)}
\end{pspicture}

\end{center}
%\end{pspicture}
\caption{\label{f2}\footnotesize {An illustrations for the proof of
Subcase 1.2a of the theorem.}}
\end{figure}

If there is a fault-free neighbor $z_L$ of $u_L$ not in $P_1'$,
without loss of generality, assume that $z_L$ is in $P_2'$, and
$P_2'$ connect $s_2$ and $t_2$ (see Fig.~\ref{f2}). Then $z_L\ne
s_2$ since $z_L\in Y$ and $s_2\in X$. Let $z_L'$ be the neighbor of
$z_L$ in $P_2'(s_2,z_L)$ (maybe $z_L'=s_2$), then $z_R'$ is
fault-free in $R$. Since both $v_L$ and $z_L'$ are in $X$, both
$v_R$ and $z_R'$ are in $Y$, and so the distance between $z_R'$ and
$v_R$ is even. Let $P_R$ be a $z_R'v_R$-path with at least
$2^{n-1}-2f_R-1$ vertices in $R$, and let
 $$
 P_1^*=s_1z_L+P_2'(z_L,t_2),\
 P_2^*=P_2'(s_2,z_L')+z_L'z_R'+P_R+v_Rv_L+P_1'(v_L,t_1).
 $$
Replacing $P_1$ and $P_2$ in (\ref{e3.2}) by $P_1^*$ and $P_2^*$
yields $k$ disjoint fault-free $ST$-paths with vertices at least
 $$
 (2^{n-1}-2(f_L-1)-|\{w\}|)+(2^{n-1}-2f_R-1)=2^n-2f.
 $$

Summing up the above discussion, the theorem holds if $f_R=1$.

{\it Subcase 1.2b.}\quad $f_R=0$.

Let $S'=S-\{s_k\}$ and $T'=T-\{t_k\}$. Since
$f_L=f\leq 2n-2k-3=2(n-1)-2(k-1)-3$, by the induction hypothesis, in
$L$ there are $(k-1)$ disjoint fault-free $S'T'$-paths
 \begin{equation}\label{e3.3}
 P_1',P_2',\ldots,P_{k-1}'
 \end{equation}
containing at least $2^{n-1}-2f$ vertices. Let
 $$
 K=V(P_1'\cup P_2'\cup\ldots\cup P_{k-1}'),
 $$
and let $s_R$ and $t_R$ be neighbors of $s_k$ and $t_k$ in $R$,
respectively. Since $s_k$ and $t_k$ are in different partite sets in
$Q_n$, $d_{L}(s_k,t_k)$ is odd, and so is $d_{R}(s_R,t_R)$.

If $|\{s_k,t_k\}\cap K|=0$ then, since $d_R(s_R,t_R)$ is odd,  by
Lemma~\ref{lem2.1}, there is an $s_Rt_R$-path $P_R$ containing
$2^{n-1}$ vertices in $R$. Let
 \begin{equation}\label{e3.4}
 \begin{array}{rl}
 & P_i=P_i',\ i=1, 2,3,\ldots,k-1,\\
 &P_k=s_ks_R+P_R+t_Rt_k.
 \end{array}
 \end{equation}
Then $P_1,\ldots,P_{k-1},P_k$ are $k$ disjoint fault-free $ST$-paths
containing at least $2^n-2f$ vertices, as required.

If $|\{s_k,t_k\}\cap K|=2$, then there are two cases to be
considered.

a)\ $s_k$ and $t_k$ are both in the same path in (\ref{e3.3}), say
$P_1'$. Let $x_L$ and $y_L$ be two neighbors of $s_k$ and $t_k$ in
$P_1'$ but not in the subpath $P_1'(s_k,t_k)$, respectively. Since
$d_L(s_k,t_k)$ is odd, both $d_L(x_L,y_L)$ and $d_R(x_R,y_R)$ are
odd. By Lemma~\ref{lem2.1}, there is an $x_Ry_R$-path $P_R$
containing $2^{n-1}$ vertices in $R$. Let
 $$
 \begin{array}{rl}
 & P_1^*=P_1'(s_1,x_L)+x_Lx_R+P_R+y_Ry_L+P_1'(y_L,t_1),\\
 & P_k^*=P_1'(s_k,t_k).
 \end{array}
 $$
Replacing $P_1$ and $P_k$ in (\ref{e3.4}) by $P_1^*$ and $P_k^*$
yields $k$ disjoint fault-free $ST$-paths containing at least
$2^n-2f$ vertices.

b)\ $s_k$ and $t_k$ are in different paths in (\ref{e3.3}). Without
loss of generality, suppose that $s_k$ is in $P_1'$ and $t_k$ is in
$P_2'$, two end-vertices of $P_i'$ are $s_i$ and $t_i$ for $i=1,2$.
Let $x_L$ be the neighbor of $s_k$ in $P_1'(s_1,s_k)$, and $y_L$ be
the neighbor of $t_k$ in $P_2'(t_k,t_2)$. Since $d_L(s_k,t_k)$ is
odd, both $d_L(x_L,y_L)$ and $d_R(x_R,y_R)$ are odd. By
Lemma~\ref{lem2.1}, there is an $x_Ry_R$-path $P_R$ containing
$2^{n-1}$ vertices in $R$. Let
 $$
 \begin{array}{rl}
 &P_1^*=P_1'(s_1,x_L)+x_Lx_R+P_R+y_Ry_L+P_2'(y_L,t_2),\\
 &P_2^*=P_2'(s_2,t_k),\\
 &P_k^*=P_1'(s_k,t_1).
 \end{array}
 $$
Replacing $P_1, P_2$ and $P_k$ in (\ref{e3.4}) by $P_1^*, P_2^*$ and
$P_k^*$ yields $k$ disjoint fault-free $ST$-paths containing at
least $2^n-2f$ vertices.

The remaining case is $|\{s_k,t_k\}\cap K|=1$. Without loss of
generality, assume $s_k\in K$ and $t_k\notin K$. We can further
assume that $s_k$ is in $P_1'$ with two end-vertices $s_1$ and
$t_1$. Let $x_L$ be the neighbor of $s_k$ in $P_1'(s_1,s_k)$, and
$t_R$ be the neighbor of $t_k$ in $R$. Since $f_R=0$ and $x_R$ and
$t_R$ both are in $X$, there is a free-fault $x_Rt_R$-path $P_R$
containing $2^{n-1}-1$. Let
 $$\begin{array}{rl}
 &P_1^*=P_1'(s_1,x_L)+x_Lx_R+P_R+t_Rt_k\\
 &P_k^*=P_1'(s_k,t_1).
 \end{array}
 $$
Replacing $P_1$ and $P_k$ in (\ref{e3.4}) by $P_1^*$ and $P_k^*$
yields $k$ disjoint fault-free $ST$-paths containing vertices at
least
 $$
 (2^{n-1}-2f+|\{t_k\}|)+(2^{n-1}-1)=2^n-2f.
 $$

Summing up the above discussion, the theorem holds when $q=k$.
\vskip6pt  For sake of convenience, we assume $f_L \le f_R$ in the
following two cases. \vskip6pt {\bf Case 2.}\ $1\le q\leq k-1$ or
$p> q =0$ and $f_R\le 2n-2k-5$.

Let $W_L$ be the set of neighbors of vertices of $S_R$ in $L$. Then
$W_L\subset Y$ since $S\subset X$. For $n\geq 5$ and $k\geq 2$, we
have that
 $$
 \begin{array}{rl}
 |Y\cap V(L)|-|T_L|-|W_L|-f &= 2^{n-2}-q-(k-p)-f\\
  &\geq 2^{n-2}-q-(k-p)-(2n-2k-3)\\
  &=(2^{n-2}-2n+k+3)+(p-q)\\
  &>p-q.
 \end{array}
 $$
This implies that there is a set $U_L$ of $(p-q)$ fault-free
vertices in $Y\cap V(L-T_L)$ such that its neighbor-set $U_R$ in $R$
is in $X\cap V(R-S_R)$ and fault-free. If $p>q$, let
\begin{equation}\label{e3.5}
U_L=\{u_1,u_2,\ldots,u_{p-q}\}\ \ {\rm and}\ \
U_R=\{v_1,v_2,\ldots,v_{p-q}\},
 \end{equation}
be any two such vertex-sets, where $u_iv_i\in E(Q_n)$ for each
$i=1,2,\ldots,p-q$. If $p=q$, let $U_L=U_R=\emptyset$.

Since $f_L\le f_R$, then $f_L \le \frac{1}{2} f \le
\frac{1}{2}(2n-2k-3)\le 2(n-1)-2k-3$. By the induction hypothesis,
in $L$ there are $p(\le k)$ disjoint fault-free paths
 $$
P_1',\ldots, P_{p-q}', P_{p-q+1},\ldots,P_p
 $$
connecting $S_L$ and $T_L\cup U_L$, and containing at least
$2^{n-1}-2f_L$ vertices, where $P_i'$ connect some vertex in $S_L$
to the vertex $u_i$ in $U_L$ for each $i=1,2,\ldots,p-q$ provided
that $p>q$.

If $q\geq 1$, then $f_R\leq f\leq 2n-2k-3\leq 2(n-1)-2(k-q)-3$ and
$k-q\leq k-1\le n-3$. If $q=0$ and $f_R\leq 2n-2k-5$, then
$f_R\leq 2(n-1)-2k-3, k\leq n-3$. In
any case, by the induction hypothesis, in $R$ there are $k-q$
disjoint fault-free paths
 $$
P_{1}'',\ldots, P_{p-q}'', P_{p+1},\ldots,P_k
 $$
connecting $S_R\cup U_R$ and $T_R$, and containing at least
$2^{n-1}-2f_R$ vertices, where $P_i''$ connect the vertex $v_i$ in
$U_R$ to some vertex in $T_R$ for each $i=1,2,\ldots,p-q$ provided
that $p>q$.
 Let
 $$
P_{i}=P_{i}'+u_iv_i+P_{i}''\ \ i=1,2,\ldots,p-q.
 $$
Then $k$ paths $P_1,P_2,\ldots,P_{p-q},P_{p-q+1},\ldots,P_p,
P_{p+1},\ldots, P_k$ satisfy our requirements.

\vskip6pt {\bf Case 3.}\  $p> q =0$ and $f_R\ge 2n-2k-4$.

In this case, $T_R=T$ and $f_L\leq 1$. We
can write (\ref{e3.5}) as
 \begin{equation}\label{e3.6}
 U_L=\{u_1,u_2,\ldots,u_p\}\ \ {\rm and}\ \
 U_R=\{v_1,v_2,\ldots,v_{p}\}.
 \end{equation}

Since $f_L\le 1$, by induction, in $L$ there are $p$
disjoint fault-free $S_LU_L$-paths
  \begin{equation}\label{e3.7}
 P_1',P_2',\ldots, P_{p}'
 \end{equation}
containing at least $2^{n-1}-2f_L$ vertices, without loss of generality, say $P_i'$ connect
$s_i$ to $u_i$ for each $i=1,2,\ldots,p$.

Since $f_R\leq 2n-2k-3= 2(n-1)-2(k-1)-3$, by the induction
hypothesis, there are $k-1$ disjoint fault-free paths
  \begin{equation}\label{e3.8}
 P''_{2},\ldots, P''_p, P_{p+1},\ldots,P_k
  \end{equation}
connecting $S_R\cup U_R-\{v_1\}$ and $T-\{t_1\}$ in $R$, and
containing at least $2^{n-1}-2f_R$ vertices, where $P''_i$ connect the
vertex $v_i$ in $U_R$ and some vertex in $T-\{t_1\}$ for
each $i=2,3,\ldots,k$. Then
 \begin{equation}\label{e3.9}
 \begin{array}{rl}
 & P_1=P_1', \\
 & P_{i}=P_{i}'+u_iv_i+P''_{i},\ \ i=2,3, \ldots,p\\
 & P_{i}, \ \ i=p+1,\ldots ,k
 \end{array}
 \end{equation}
are $k$ disjoint fault-free paths between $S$ and
$(T-\{t_1\})\cup\{u_1\}$ containing at least $2^n-2f$
vertices in $Q_n$. Without loss of generality, say $P_i$ connect
$s_i$ to $t_i$ for each $i=2,3,\ldots,k$, $P_1$ connect $s_1$
to $u_1$. Let
 $$
 B=V(P_{2}\cup\ldots\cup P_k).
 $$
 Note if $p=1$, $P'_i$ in (\ref{e3.7})  and $P''_i$ in (\ref{e3.8})
are empty, and we can choose $t_1$ such that $s_1t_1\not\in E(Q_n)$.

There are two subcases to be considered.

a)\ $t_1\notin B$.  Let $t_L$ be the neighbor of $t_1$ in $L$.

Suppose that $t_L$ is fault-free. If  $t_L\in S $, say $t_L=s_1$.
By assumption $s_1t_1\not\in E(Q_n)$ if $p=1$, then $p \ge 2$. We consider
 $t_L$ a temporary faulty vertex. Let $S'=S-\{s_1\},
U'=U_L-\{u_1\}$. Then $L$ contains at most two faulty
vertices, $p-1$ disjoint fault-free $S'U'$-paths $P'_2,\ldots,P'_p$
 shown in (\ref{e3.7}) contain at least $2^{n-1}-2(f_L+1)$ vertices,
 $P_2,P_3,\ldots,P_k$ in  (\ref{e3.9}) are $k-1$ disjoint paths containing
 at least $2^n-2f-2$ vertices. Together with the path $ P_1=s_1t_1$, obtain
 the $k$ paths satisfied the requirements.
If $t_L\not \in S $, since $f_L\le1$,
$|F_L\cup (U_L- \{u_1\})|\le k \le n-3$,  thus $t_L$ has at least
a fault-free neighbor $z\in V(L)-(U_L- \{u_1\})$. Since $z\in
Y$, we can choose $u_1$ of $U_L$ be $z$ in (\ref{e3.6}). We consider $t_L$
a temporary faulty vertex. Then $L$ contains at most two faulty
vertices, $p$ disjoint fault-free $S_LU_L$-paths shown in (\ref{e3.7})
contain at least $2^{n-1}-2(f_L+1)$ vertices. Let
 $$
 P_1^*=P_1+u_1t_L+t_Lt_1.
 $$
Replacing $P_1$ in (\ref{e3.9}) by $P_1^*$ yields $k$ disjoint
fault-free $ST$-paths in $Q_n$. Note that $t_1$ and $t_L$ are two
new vertices. Then $k$ paths $P_1^*,P_2,\ldots,P_k$ satisfy our
requirements.

Now suppose that $t_L$ is a faulty vertex. Then $t_L$ is the only
faulty vertex in $L$. Take a fault-free neighbor, say $w_R$ of $t_1$
in $R$. Then its neighbor $w_L$ in $L$ is fault-free and in $Y$. If
$w_R\notin B$, then choose $u_1=w_L\in U_L$ in (\ref{e3.6}). Let
 $$
 P_1^*=P_1+u_1w_R+w_Rt_1.
 $$
Replacing $P_1$ in (\ref{e3.9}) by $P_1^*$ yields $k$ disjoint
fault-free $ST$-paths in $Q_n$, as required. Assume that $w_R$ is in
some $P_i$ $(2\leq i\leq k)$, say in $P_2(s_2,t_2)$. Then $w_R \ne t_2$ since
$w_R$ and $t_2$ are in different partite sets. Let $u_R$ be the
neighbor of $w_R$ in $P_2(w_R,t_2)$ and $u_L$ be the neighbor
of $u_R$ in $L$. Then $u_L\in X$.

 If $u_L\in S$, say $u_L=s_1$. If $p=1$, we can choose
 another $w_R$ such that $u_L\not=s_1$, so let $p\ge2$.
 We consider $s_1$ be a temporary fault vertex. Let $S'=S-\{s_1\},
U'=U_L-\{u_1\}$. Then $L$ contains at most two faulty
vertices, there are $p-1$ disjoint fault-free $S'U'$-paths $P'_2,\ldots,P'_p$
 containing at least $2^{n-1}-2(f_L+1)$ vertices,
 $P_2,P_3,\ldots,P_k$ in  (\ref{e3.9}) are $k-1$ disjoint paths containing
 at least $2^n-2f-2$ vertices.
 Let  $$
\begin{array}{rl}
& P_1=s_1u_R+P_2(u_R,t_2),\\
& P_2^*=P_2(s_2,w_R)+w_Rt_1.
 \end{array}$$
 Thus  $P_1,P^*_2,P_3,\ldots,P_k$ are $k$ disjoint fault-free $ST$-paths in $Q_n$,
where $t_1$ and $s_1$ are two new vertices, as required.

If $u_L \not \in S$, choose a fault-free neighbor $z_L$
of $u_L$ in $L$ such that it is not any $u_i$ in $U_L$ for $2\leq
i\leq p$. Then $z_L\in Y$. So choose $u_1=z_L\in U_L$ in
(\ref{e3.6}). We consider $u_L$ a temporary faulty vertex. Then $L$
contains at most two faulty vertices. $k$ disjoint fault-free
$SU_L$-paths shown in (\ref{e3.7}) contain at least
$2^{n-1}-2(f_L+1)$ vertices. Let
 $$
 \begin{array}{rl}
 &P_1^*=P_1+z_Lu_L+u_Lu_R+P_2(u_R,t_2),\\
 &P_2^*=P_2(s_2,w_R)+w_Rt_1.
 \end{array}
 $$
Replacing $P_1$ and $P_2$ in (\ref{e3.9}) by $P_1^*$ and
$P_2^*$ yields $k$ disjoint fault-free $ST$-paths in $Q_n$,
where $t_1$ and $u_L$ are two new vertices, as required.

b)\ $t_1\in B$.

In this case, without loss of generality, we assume $t_1$ is in
$P_2$ that connects $s_2$ to $t_2$. Let $u_R$ and $v_R$ be two
neighbors of $t_1$ in $P_2$ , and let $u_L$ and $v_L$ be two
neighbors of $u_R$ and $v_R$ in $L$, respectively. Then $u_L$ and
$v_L$ both are in $Y$. Without loss of generality, assume that $v_R$
is in $P_2(s_2,t_1)$ and $u_R$ is in $P_2(t_1, t_2)$.

If $f_L = 0$ or $u_L$ is fault-free, then choose $u_1=u_L$ in $U_L$
in (\ref{e3.6}). Let
 $$
 \begin{array}{rl}
 & P^*_1 = P_2(s_2, t_1),\\
 & P^*_2 = P_1 + u_Lu_R + P_2(u_R, t_2).
 \end{array}
 $$
Replacing $P_1$ and $P_2$ in (\ref{e3.9}) by $P^*_1$ and $P^*_2$
yields $k$ disjoint fault-free $ST$-paths in $Q_n$, as required.

Assume now that $f_L = 1$ and $u_L$ is fault vertex. Since $f\le
2n-2k-3$, $f_R=2n-2k-4$. If $d(u_R,t_2)>1$, we can choose $t_2$
instead of $t_1$, which comes down to the above case, yields $k$
disjoint fault-free $ST$-paths as required. Thus, we can assume
$d(u_R,t_2)=1$.

Let $t_{2L}$ be the neighbor of $t_2$ in $L$, then $u_Lt_{2L}\in E(Q_n)$.
If $t_{2L}\in S_L$, say $t_{2L}=s_1$. By assumption $s_1t_1\not\in E(Q_n)$ if $p=1$,
 then $p>1$. Let $S'=S-\{s_1\},
U'=U_L-\{u_1\}$.  Using Lemma~\ref{lem2.4}
and choosing $xy$ to be $u_Lt_{2L}$, then there is $p-1\, (\le k)$
disjoint $S'U'$-paths
  \begin{equation}\label{e3.10}
  P_2',\ldots, P_{p}'
 \end{equation}
that contain all vertices in $L-\{u_L,t_{2L}\}$, say  $P_i'$
connect $s_i$ to $u_i$ for each $i=2,3,\ldots,p$.
Combining (\ref{e3.8}) with (\ref{e3.10}), we
obtain $k-1$ disjoint fault-free paths between $S'$ and
$T-\{t_1\}$
 \begin{equation}\label{e3.11}
 \begin{array}{rl}
 & P_{i}=P_{i}'+u_iv_i+P''_{i},\ \ i=2,3,\ldots,p\\
 & P_{i}\ \ i=p+1,\ldots,k
 \end{array}
 \end{equation}
that contain at least
 $
 (2^{n-1}-2)+(2^{n-1}-2f_R)=2^n-2f
 $
vertices in $Q_n$. Let
 $$
 P_1=s_1t_2,\ \ {\rm and }\ \ P^*_2=P_2(s_2,t_1).
 $$
then the paths $P_1,P^*_2,P_3,\ldots,P_k$ as required.

If  $t_{2L}\not \in S_L$, let $z \in N(t_{2L})\cap L- \{u_L \}-U'$.
Choose $u_1=z$ in $U_L$ in (\ref{e3.6}).
 Using Lemma~\ref{lem2.4}
and choosing $xy$ to be $u_Lt_{2L}$, then there is $p\, (\le k)$
disjoint $S_LU_L$-paths
  \begin{equation}\label{e3.12}
 P_1',P_2',\ldots, P_{p}'
 \end{equation}
that contain all vertices in $L-\{u_L,t_{2L}\}$, say  $P_i'$
connect $s_i$ to $u_i$ for each
$i=1,2,\ldots,p$. Combining (\ref{e3.8}) with (\ref{e3.12}), we
obtain $k$ disjoint fault-free paths between $S$ and
$(T-\{t_1\})\cup\{u_1\}$
 \begin{equation}\label{e3.13}
 \begin{array}{rl}
 & P_1=P_1', \\
 & P_{i}=P_{i}'+u_iv_i+P''_{i},\ \ i=2,3,\ldots,p\\
 & P_{i},\ \ i=p+1,\ldots,k\ ({\rm defined\ in}\ (\ref{e3.8}))
 \end{array}
 \end{equation}
that contain at least
 $
 (2^{n-1}-2)+(2^{n-1}-2f_R)=2^n-2f
 $
vertices in $Q_n$. Let
 $$
 P^*_1=P_1+u_1t_{2L}+t_{2L}t_2\ \ {\rm and }\ \ P^*_2=P_2(s_2,t_1).
 $$
Replacing $P_1$ and $P_2$ in (\ref{e3.13}) by $P^*_1$ and $P^*_2$
yields $k$ disjoint fault-free $ST$-paths in $Q_n$. We remove one
vertex $u_R$ from paths in (\ref{e3.11}), and add one vertex $t_{2L}$
to obtain new paths. Thus, these
paths still contain at least $2^n-2f$ vertices in $Q_n$, as
required.

The theorem follows.

%%%%%%%%%%%%%%%%%%%%%%%%%%%%%%%%%%%%%%%%%%%%%%%%%%%%%%%%%%%%%%%%%%%%%%%%%%%%%%%%%%

\end{document}